\documentclass{article}
\usepackage{amsmath,amssymb}
\usepackage{theorem}
\usepackage{color}

\theorembodyfont{\itshape}
\theoremstyle{plain}
\newtheorem{theorem}{Theorem}[section]
\newtheorem{lemma}{Lemma}[section]
\newtheorem{corollary}{Corollary}[section]
\newtheorem{proposition}{Proposition}[section]
\theorembodyfont{\rmfamily}

\newtheorem{question}{Question}[section]
\newtheorem{problem}{Problem}[section]

\newtheorem{remark}{Remark}[section]

\newcommand{\Cl}{\mbox{{\rm Cl}}}

\begin{document}

\title{\bf On $\sigma$-compact Hattori spaces.}

\author{Vitalij A.~Chatyrko}


\maketitle

\begin{abstract} We present several characterizations of $\sigma$-compact Hattori spaces, and 
reject some  possible characterization candidates  of the spaces.
\end{abstract} 


\smallskip
{\it 2000 AMS (MOS) Subj. Class.:} Primary 54A10   
\medskip
\baselineskip=18pt

\section{Introduction}

Let $\mathbb R$ be the set of real numbers and $A$ be a subset of $\mathbb R$.

In \cite{H} Hattori introduced a topology $\tau(A)$ on $\mathbb R$ defined as follows:
\begin{itemize}
	\item[(1)] if $x \in A$ then $\{(x-\epsilon, x+\epsilon): \epsilon > 0\}$ is a nbd open basis at $x$, 
	\item[(2)] if $x \in \mathbb R \setminus A$ then $\{[x, x+\epsilon): \epsilon > 0\}$ is a nbd open basis at $x$.
\end{itemize}

Note that $\tau(\emptyset)$ (respectively, $\tau(\mathbb R)$) is the Sorgenfrey topology $\tau_S$ (respectively, the Euclidean topology $\tau_E$) on the reals.

The topological spaces $(\mathbb R, \tau(A)), A \subseteq \mathbb R,$ are called {\it Hattori spaces} and denoted by $H(A)$ or $H$ (if $A$ is unimportant for a discussion). It is easy to see that the identity mapping of reals is a continuous bijection of any $H$-space onto the real line.

Let us recall (\cite{CH}) that every $H$-space is $T_1$, regular, hereditary Lindel\"of and hereditary separable. However there are topological properties as the metrizability or the Cech-completeness which some $H$-spaces  possess and other $H$-spaces do not possess. When the $H$-spaces possess these properties one can find in \cite{K} and \cite{BS}. 

Recall (\cite{EJ}) that each compact subset of the Sorgenfrey line $H(\emptyset)$ is countable. So 
the space $H(\emptyset)$  cannot be $\sigma$-compact unlike to the space $H(\mathbb R)$  
(the real line) which is evidently $\sigma$-compact. 

The following natural question was posed by F. Lin and J. Li.

\begin{question}(\cite[Question 3.7]{LL}) For what subsets $A$ of $\mathbb R$  are the spaces $H(A)$ $\sigma$-compact?
\end{question}

F. Lin and J. Li  also noted
 
\begin{proposition}\label{prop_1}(\cite[Theorem 3.13]{LL}) For an arbitrary subset $A$ of $\mathbb R$, if $H(A)$ is $\sigma$-compact, then
$\mathbb R \setminus A$ is countable and nowhere dense in $H(A)$.
 $\Box$ 
\end{proposition}

\begin{proposition}\label{prop_2}(\cite[Theorem 3.14]{LL}) For an arbitrary subset $A$ of $\mathbb R$, if $\mathbb R \setminus A$ is countable and scattered
in $H(A)$, then $H(A)$ is $\sigma$-compact.
 $\Box$ 
\end{proposition}

In this note I present several characterizations of Hattori spaces to be $\sigma$-compact, and show that 
the implications of Propositions \ref{prop_1} and \ref{prop_2} are not invertible. Moreover, Proposition \ref{prop_2} (formulated as above) does not hold, its corrected version is presented in Corollary \ref{cor}.

For standard notions we refer to \cite{E}.

\section{Main results}

First of all let us recall the following fact.
\begin{lemma}\label{lem_CH}(\cite[Lemma 2.1]{CH}) Let $A \subseteq \mathbb R$ and $B \subseteq A$ and $C \subseteq \mathbb R \setminus A$. Then 
\begin{itemize}
	\item[(i)] $\tau(A)|_B = \tau_E|_B$, where $\tau_E$ is the Euclidean topology on $\mathbb R$, and 
	\item[(ii)] $\tau(A)|_C = \tau_S|_C$, where $\tau_S$ is the Sorgenfrey topology on $\mathbb R$.  $\Box$ 
\end{itemize}

\end{lemma}

\begin{proposition}\label{prop_3} For an arbitrary subset $A$ of $\mathbb R$, if $B= \mathbb R \setminus A$ is countable and it is a $G_\delta$-subset of
the real line (in particular, if $B$ is countable and closed in the real line), then $H(A)$ is $\sigma$-compact.
 
\end{proposition}

Proof.  Let us note that on the real line our set $A$ is an $F_\sigma$-set and hence it is $\sigma$-compact there.
So by Lemma \ref{lem_CH} $A$ is  $\sigma$-compact in $H(A)$ too. Since $B$ is countable we get that $H(A)$ is $\sigma$-compact. $\Box$

Since every scattered subset of the real line is a $G_\delta$ (see \cite[Corollary 4]{KR}) we get the following.

\begin{corollary}\label{cor} For any subset $A$ of $\mathbb R$, if $\mathbb R \setminus A$ is countable and 
scattered in
the real line, then $H(A)$ is $\sigma$-compact.  $\Box$ 
\end{corollary}

We continue with several characterizations of $H$-spaces to be $\sigma$-compact.
 
 \begin{theorem}\label{character} Let $A \subseteq \mathbb R$ and $B = \mathbb R \setminus A$. Then
the following conditions are equivalent.
\begin{itemize}
\item[(a)] There exist a $\sigma$-compact subset $D$ and a closed subset $C$ of the space  $H(A)$ 
such that $B \subseteq C \subseteq D$.
\item[(b)] There exists a closed $\sigma$-compact subset $C$ of the space $H(A)$ 
such that $B \subseteq C$.
\item[(c)] The closure $\Cl_{H(\mathbb R)}(B)$ of $B$ in the real line is 
$\sigma$-compact in $H(A)$.
\item[(d)] The closure $\Cl_{H(A)}(B)$ of $B$ in the space $H(A)$ is 
$\sigma$-compact in $H(A)$.
\item[(e)] the space $H(A)$ is $\sigma$-compact.

\end{itemize}

 \end{theorem}
 
 Proof. The following implications are obvious: $(e) => (a)$, $(a) => (b)$, $(c) => (b)$, 
 $(b) => (d)$,  $(e) => (c)$.

 Let us show $(d) => (e)$. Since $B \subseteq \Cl_{H(A)}(B)$
each point $x \in H(A) \setminus \Cl_{H(A)}(B)$ has inside the set $H(A) \setminus \Cl_{H(A)}(B)$ an open nbd which is  an open interval of the  real line. Since the space $H(A)$ is hereditarily Lindel\"of  the set $H(A) \setminus \Cl_{H(A)}$ is a $\sigma$-compact subset of $H(A)$  (see Lemma \ref{lem_CH}).
Thus even $H(A)$ is $\sigma$-compact. $\Box$

\begin{remark} Note that the set  $\Cl_{H(\mathbb R)}(B)$ does not need to be $\sigma$-compact in the space $H(A)$ (it is of course closed there) even if it is compact in the real line, see  Proposition \ref{my_prop_2}.
\end{remark}

Let us consider in the set of reals the standard Cantor set $\mathbb C$ on the closed interval $[0,1]$ which can be defined as follows.

For any closed bounded  interval $[a, b]$ of $\mathbb R$  put $$F([a,b]) = \{[a, \frac{2}{3}a + \frac{1}{3}b], [\frac{1}{3}a + \frac{2}{3}b,b]\}.$$
Then for each $n \geq 0$ by induction define  a family $\mathcal C_n$ of closed intervals:  
$$\mathcal C_0 = \{[0,1]\},  \ \mathcal C_n = \{F([a,b]): [a,b] \in \mathcal C_{n-1}\}.$$ 
The standard Cantor set $\mathbb C$ of the closed interval $[0,1]$ is the intersection $\cap_{n=0}^\infty (\cup \mathcal C_n)$, where $\cup \mathcal C_n$ is the union of all closed intervals from the family $\mathcal C_n$.

Put now $B_1 = \{a: [a,b] \in \mathcal C_n, n \geq 0\}$, $B_2 = \{b: [a,b] \in \mathcal C_n, n \geq 0\}$ and $A_1 = \mathbb R \setminus B_1$, $A_2 = \mathbb R \setminus B_2$. We will use the notations below.

\begin{remark} Let us note that on the real line (i.e. on the reals with  the Euclidean topology) the set $\mathbb C$ is compact, the sets $B_1$ and $B_2$ (which are subsets of $\mathbb C$) are homeomorphic to the space of rational numbers $\mathbb Q$, the sets 
$\mathbb C \setminus B_1$ and $\mathbb C \setminus B_2$ are homeomorphic to the space of irrational numbers $\mathbb P$. Moreover,  $B_1$ and $B_2$ are nowhere dense in the real line. 
\end{remark}

\begin{remark}\label{rem_2} Let us note that a set $Y \subset \mathbb R$ is nowhere dense in the real line iff $Y$ is nowhere dense in any $H$-space 
(see for example, \cite[Lemma 3.3]{CN}).
\end{remark}

\begin{proposition} For the space $H(A_1)$ the following is valid.
\begin{itemize}
\item[(a)] The subspace $B_1$ of $H(A_1)$ is nowhere dense in $H(A_1)$ and it is  homeomorphic to the space of rational numbers $\mathbb Q.$
\item[(b)] The subspace $\Cl_{H(A_1)}(B_1)$ of $H(A_1)$  is homeomorphic to the standard Cantor set $\mathbb C$ on the real line, and the subspace $\Cl_{H(A_1)}(B_1) \setminus B_1$ of $H(A_1)$  is homeomorphic to the space of irrational numbers $\mathbb P$.  
\item[(c)] The space $H(A_1)$ is $\sigma$-compact.
\end{itemize}
\end{proposition}

Proof. (a) and (b) are obvious. Theorem \ref{character} and (b) prove (c). $\Box$

\begin{corollary} Proposition \ref{prop_2} is not invertible. $\Box$
\end{corollary}

Proof. Let us note that $H(A_1)$ is $\sigma$-compact but the subspace $B_1$ of $H(A_1)$ is not scattered. $\Box$

\begin{corollary} Proposition \ref{prop_3} is not invertible. 
\end{corollary}

Proof. Let us note that $H(A_1)$ is $\sigma$-compact but $B_1$ is not a $G_\delta$-subset of the Cantor set $\mathbb C$ in the real line and hence it is not a $G_\delta$ in the real line. $\Box$

\begin{proposition}\label{my_prop_2} For the space $H(A_2)$ the following is valid.
\begin{itemize}
\item[(a)] The subspace $B_2$ of $H(A_2)$ is nowhere dense in $H(A_2)$ and it is  homeomorphic to the space of natural numbers $\mathbb N.$
\item[(b)] The subspace $\Cl_{H(A_2)}(B_2)$ of $H(A_2)$  is equal to the standard Cantor set $\mathbb C$ of $\mathbb R$,  and it is not $\sigma$-compact. The subspace $\Cl_{H(A_2)}(B_2) \setminus B_2$ of $H(A_2)$  is homeomorphic to the space of irrational numbers $\mathbb P$, 
\item[(c)] The space $H(A_2)$ is not $\sigma$-compact.
\end{itemize}
\end{proposition}

Proof. (a) is obvious. 

In (b) let us show that  the subspace $\Cl_{H(A_2)}(B_2)$ of $H(A_2)$  is not $\sigma$-compact. Assume that the subspace $\Cl_{H(A_2)}(B_2)$ of $H(A_2)$ is $\sigma$-compact, i.e. $\Cl_{H(A_2)}(B_2) = \cup_{i=1}^\infty K_i$, where $K_i$ is compact in $H(A_2)$. Note that for each $i$ the set $K_i$ is compact in the real line and the Cantor set $\mathbb C$ with the topology from the real line is the union $\cup_{i=1}^\infty K_i$. Hence there is an open interval $(c, d)$ of the reals and some $i$ such that 
$(c, d) \cap \mathbb C \subseteq K_i$. Moreover, there exist  points $b_0, b_1, \dots$ of $B_2$ such that 
$b_1 < b_2 < \dots < b_0$ and the sequence  $\{b_j\}_{j=1}^\infty$ tends to $b_0$ in the real line.
Since at the points of $B_2$ the topology of $H(A_2)$ is  the Sorgenfrey topology  we get a contradiction with the compactness of $K_i$ in the space $H(A_2)$.  

Theorem \ref{character} and (b) prove (c). $\Box$

\begin{corollary} Proposition \ref{prop_1} is not invertible. 

\end{corollary}

Proof. Let us note that $B_2$ is nowhere dense in $H(A_2)$ (see Remark \ref{rem_2}) but the space 
$H(A_2)$ is not $\sigma$-compact. $\Box$

\begin{corollary} Proposition \ref{prop_2} does not hold. $\Box$

Proof. Let us note that $B_2$ is scattered in $H(A_2)$ and the space $H(A_2)$ is not  $\sigma$-compact. 
$\Box$

\end{corollary}

\section{Additional questions}

The following is obvious.
\begin{itemize}
\item[(a)] If a space $X$ is $\sigma$-compact then a subset $Y$ of $X$ is $\sigma$-compact iff it is an $F_\sigma$-subset of $X$. In particular, a subset of the real line is $\sigma$-compact iff it is an $F_\sigma$-set.
\item[(b)] A subset of the Sorgenfrey line is $\sigma$-compact iff it is countable, 
\item[(c)]  A subset of the space $\mathbb P$ of irrational numbers is $\sigma$-compact iff
it is homeomorphic to an $F_\sigma$-subset of the standard Cantor set $\mathbb C$ on the real line.

\end{itemize}

One can pose the following problem.
\begin{problem}\label{problem} Let $A \subseteq \mathbb R$. Describe the $\sigma$-compact subsets of $H(A)$.
\end{problem}

Let us note in advance that according to (a) if $H(A)$ is $\sigma$-compact then a subset of $H(A)$ is $\sigma$-compact iff
it is an $F_\sigma$-subset of $H(A)$ 

Below we present some other answers to Problem \ref{problem} by the use of observations (b) and (c) and some known facts.

\begin{proposition}(\cite[Theorem 6]{K} and \cite[Theorem 2.8]{BS}) $H(A)$ is homeomorphic to the Sorgenfrey line iff $A$ is scattered. $\Box$
\end{proposition}

\begin{corollary} If $A$ is scattered then  a subset of $H(A)$ is $\sigma$-compact iff it is countable. $\Box$ 
\end{corollary}

\begin{proposition}(\cite[Proposition 3.6]{BS}) $H(A)$ is homeomorphic to the space $\mathbb P$ of irrational numbers iff $\mathbb R \setminus A$ is dense in the real line and  countable. $\Box$
\end{proposition}

\begin{corollary} If $\mathbb R \setminus A$ is dense in the real line and  countable then  a subset of $H(A)$ is $\sigma$-compact iff it is homeomorphic to an $F_\sigma$-subset of the standard Cantor set $\mathbb C$ on the real line. $\Box$
\end{corollary}

Since the space $H(A_2)$ from Proposition \ref{my_prop_2} is not $\sigma$-compact (as well as any subset of
$H(A_2)$ containing some $[a,b]$ from $\mathcal C_n, n = 0, 1, 2, \dots$) one can pose the following question.

\begin{question} What subsets of $H(A_2)$ are $\sigma$-compact?
\end{question}

\noindent(V.A. Chatyrko)\\
Department of Mathematics, Linkoping University, 581 83 Linkoping, Sweden.\\
vitalij.tjatyrko@liu.se


\begin{thebibliography}{99} 

\bibitem[BS]{BS}  A.~Bouziad,  E.~Sukhacheva, On Hattori spaces, Comment. Math. Univ. Carolin. 58,2 (2017) 213-223
\vskip0.3cm



\bibitem[CH]{CH} V.~A.~Chatyrko, Y.~Hattori, A poset of topologies on the set of real numbers, Comment. Math. Univ. Carolin.  54, 2 (2013)  189-196.
\vskip0.3cm

\bibitem[CN]{CN} V.~A.~Chatyrko, V.~Nyagahakwa, 
Sets with the Baire property in topologies formed from a given topology and ideals of sets, 
Questions and answers in General Topology, 35 (2017) 59-76
\vskip0.3cm



\bibitem[E]{E} R.~Engelking, General Topology, Heldermann Verlag, Berlin, 1989.
\vskip0.3cm

\bibitem[EJ]{EJ} M.~S.~Espelie, J.~E.~Joseph, Compact subspaces of the Sorgenfrey line, Math. Magazine, 49 (1976) 250-251.


\bibitem[H]{H} Y.~Hattori, Order and topological structures of poset of the formal balls on metric spaces,  Mem.Fac.Sci.Eng. Shimane Univ. Ser. B Math.  Sci., 43 (2010)  13 - 26.


\bibitem[KR]{KR} V.~Kannan, M.~Rajagopalan, On scattered spaces, Proc. Amer. Math. Soc., 43(2)(1974)
402-408.

\bibitem[K]{K} J. Kulesza, Results on spaces between the Sorgenfrey topology and the usual topology on $\mathbb R$,  
Topol. Appl. 231 (2017) 266-275
\vskip0.3cm






\bibitem[LL]{LL} F.~Lin, J.~Li, Some topological properties of spaces between the Sorgenfrey  and  usual topologies on the real numbers. arXiv: 1807.06938v4 [math.GN] 21Dec 2022











\end{thebibliography}
\end{document}